# About epimorphisms of a fundamental group of a splitting surface.

**G.S. Makanin**

J. Stollings [1] and V. Jacko [2] proved, that the following two statements are equivalent.

**I.** The closed oriented simply connected thrivariate manifold, whose Heegaard splitting has the family $g > 1$, is homomorphic to a thrivariate sphere.

**II.** If $F_g'$, $F_g''$ are free groups of the rank $g > 1$, then any epimorphism $\eta$ of a fundamental group of the splitting surface $S_g$ on the group $F_g' \times F_g''$ *passes essentially* through a free composition.

The last means that for epimorphism $\eta$ of the group $S_g$ on the group $F_g' \times F_g''$ can be found so nontrivial groups $G_1, G_2$, epimorphism $\tau$ of the group $S_g$ on the free composition $G_1 * G_2$, isomorphism $\theta$ of the group $G_1 * G_2$ on group $F_g' \times F_g''$ that, the diagram

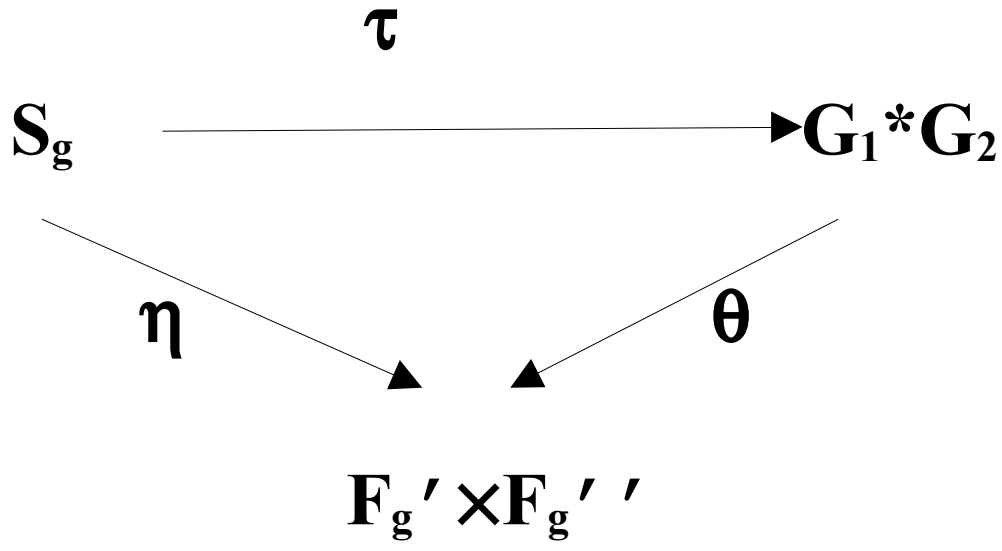

is commutative.

In this work it is proved that for every $g > 1$ any epimorphism of a fundamental group of the splitting surface $S_g$ on direct composition of two free groups of the rank $g$ passes essentially through a free composition of the group $S_{g-1}$ and Abelian group $A_2$ of the rank **2**. Thereby the hypothesis Poincaré[1], that a closed oriented simply connected thrivariate manifold is topologically equivalent to a trivariate sphere, is proved.

1 The short topological method of the proof of the hypothesis Poincaré was explained by G. Perelman

We will set the fundamental group $S_g$ of the splitting surface by generators $y_1$, $zy_2,\ldots, y_{2g-1}, y_{2g}$ and by defining relation $[y_1, y_2]\ldots[y_{2g-1}, y_{2g}]=1$

The group $S_g$ we will call $\beta$-group.

The elementary isomorphisms are determined by transformations

$y_i^{\varepsilon} \to y_i^{\varepsilon} y_j^{-\delta}$, $y_i^{\varepsilon} \to y_j^{-\delta} y_i^{\varepsilon}$, where $i,j=1,\ldots,2g$, $i \neq j$, $\varepsilon,\delta = \pm 1$.

The result of applicability of the elementary isomorphism to $\beta$-group we will call $\beta$-group.

Generators $y_1,y_2,\ldots,y_{2g-1},y_{2g}$ of the group $S_g$ have the second definition (respectively) with the letters $l_1, m_1,\ldots,l_g, m_g$.

*The Abelian group* $A_2$ of the rank $2$ is set by generators $l_1, m_1$ and defining relation $[l_1, m_1]=1$.

$F_g'$ is a free group with free generators $c_1,\ldots,c_g$, $F_g''$ is a free group with free generators $d_1,\ldots,d_g$, $F_g' \times F_g'$ is a direct composition of groups $F_g'$ and $F_g''$.

Epimorphism $\eta$ of the group $S_g$ on the group $F_g' \times F_g''$ defines epimorphism $\mu$ of the group $S_g$ on $F_g'$ and epimorphism $\nu$ of the group $S_g$ on $F_g''$ and is indicated by $\mu \times \nu$.

It is convenient to write $\beta$-group $S_g$ with epimorphism $\mu\times\nu$ in three lines

$$\begin{array}{cccccccccc}
l_1^{-1} & m_1^{-1} & l_1 & m_1 \ldots & l_g^{-1} & m_g^{-1} & l_g & m_g =1 & & S_g \\
P_1^{-1} & Q_1^{-1} & P_1 & Q_1 & P_g^{-1} & Q_g^{-1} & P_g & Q_g & & \mu \\
L_1^{-1} & M_1^{-1} & L_1 & M_1 & L_g^{-1} & M_g^{-1} & L_g & M_g & & \nu
\end{array}$$

where components $P_1,Q_1,\ldots, P_g,Q_g$ generate the group $F_g'$, and components $L_1,M_1,\ldots, L_g,M_g$ generate the group $F_g''$. Epimorphism $\mu\times\nu$ maps

$l_i \to P_i \times L_i$, $m_i \to Q_i \times M_i$ $(i=1,\ldots,g)$.

$\beta$-group $S_g$, $\mu\times\nu$ divides into composition $g$ of *commutators*

$$\begin{array}{cccc}
\ldots l_i^{-1} & m_i^{-1} & l_i & m_i \ldots \\
P_i^{-1} & Q_i^{-1} & P_i & Q_i \\
L_i^{-1} & M_i^{-1} & L_i & M_i
\end{array}$$

Let's define *commutator with the serial number* ***i*** like $[i]$.

Any epimorphism $\mu_1\times\nu_1$ of the group $S_g$ on the group $F_g' \times F_g''$, generated by the mapping of the form

$$\begin{cases}
l_1 \to c_1 \times L_1 \\
m_1 \to 1 \times M_1 \\
\ldots\ldots\ldots\ldots\ldots\ldots \\
l_g \to c_g \times L_g \\
m_g \to 1 \times M_g
\end{cases}$$

where $c_1,\ldots, c_g$ generate the group $F_g'$, $L_1,M_1 \ldots L_g,M_g$ generate the group $F_g''$, is called *exact*.

The exact epimorphism $\mu_N\times \nu_N$ of the group $S_g$ on the group $F_g' \times F_g''$ generated by the mapping of the form

$$\begin{cases} \mathbf{l_1} \rightarrow \mathbf{c_1} \times \mathbf{L_1} \\ \mathbf{m_1} \rightarrow \mathbf{1} \times \mathbf{L_1 M_1 L_1^{-1}} \\ \dots\dots\dots\dots\dots\dots \\ \mathbf{l_g} \rightarrow \mathbf{c_g} \times \mathbf{L_g} \\ \mathbf{m_g} \rightarrow \mathbf{1} \times \mathbf{L_g M_g L_g^{-1}} \end{cases}$$

is called *normal.*

The exact epimorphism $\mu_0 \times \nu_0$ of the group $\mathbf{S_g}$ on the group $\mathbf{F_g' \times F_g''}$ , generated by mapping of the form

$$\begin{cases} \mathbf{l_1} \rightarrow \mathbf{c_1} \times \mathbf{L_1} \\ \mathbf{m_1} \rightarrow \mathbf{1} \times \mathbf{M_1} \\ \dots\dots\dots\dots\dots\dots \\ \mathbf{l_g} \rightarrow \mathbf{c_g} \times \mathbf{L_g} \\ \mathbf{m_g} \rightarrow \mathbf{1} \times \mathbf{M_g} \end{cases}$$

where for every $\mathbf{i=1,...,g}$ $\mathbf{L_i\, M_i = M_i\, L_i}$ in the group $\mathbf{F_g''}$, is called *unit.*

$\gamma$-*group* is called $\beta$-group with the exact epimorphism $\mu_1 \times \nu_1$ on group $\mathbf{F_g' \times F_g''}$ .

Elementary isomorphisms in the form of

$$\mathbf{l_i^{\varepsilon} \rightarrow l_i^{\varepsilon}\, m_j^{-\delta}}, \quad \mathbf{l_i^{\varepsilon} \rightarrow m_j^{-\delta} l_i^{\varepsilon}},$$
$$\mathbf{l_i^{\varepsilon} \rightarrow l_i^{\varepsilon}\, l_j^{-\delta} \,\&\, c_i^{\varepsilon} \rightarrow c_i^{\varepsilon}\, c_j^{-\delta}}, \quad \mathbf{l_i^{\varepsilon} \rightarrow l_j^{-\delta} l_i^{\varepsilon} \,\&\, c_i^{\varepsilon} \rightarrow c_j^{-\delta} c_i^{\varepsilon}},$$
$$\mathbf{m_i^{\varepsilon} \rightarrow m_i^{\varepsilon}\, m_j^{-\delta}}, \quad \mathbf{m_i^{\varepsilon} \rightarrow m_j^{-\delta} m_i^{\varepsilon}},$$

where $\varepsilon, \delta = \pm 1$, are called $\gamma$-isomorphisms.

Besides,

*a composition of arbitrary elementary isomorphisms, which converses* $\gamma$-*group in* $\gamma$-*group, we will also call* $\gamma$-*isomorphism.*

**Lemma 1.** For group $\mathbf{S_g}$ with epimorphism $\mu \times \nu$ on group $\mathbf{F_g' \times F_g''}$ we can indicate such epimorphism $\omega$ of the group $\mathbf{S_g}$ on $\gamma$-group $\mathbf{S_g}$, such exact epimorphism $\mu_1 \times \nu_1$ of the $\gamma$-group $\mathbf{S_g}$ on group $\mathbf{F_g' \times F_g''}$ that diagram

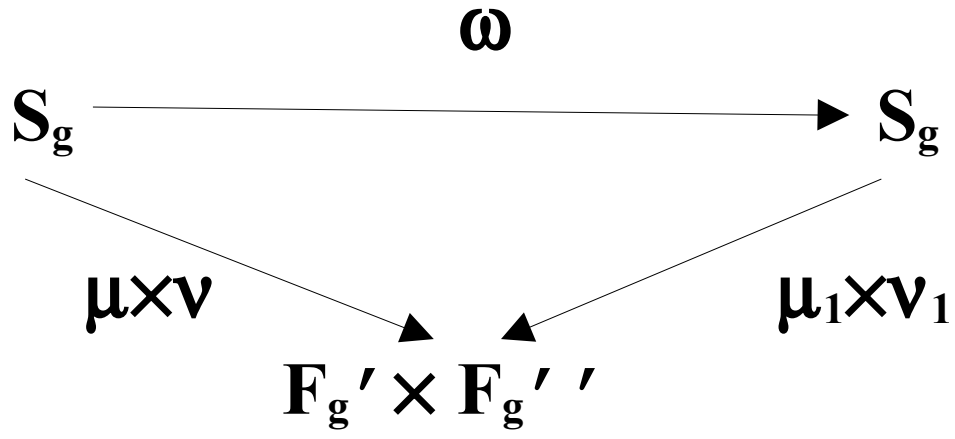

is commutative.

**Proof** * consists of conversion by elementary isomorphisms two lines

$$l_1^{-1}\quad m_1^{-1}\quad l_1\ m_1 \dots l_q^{-1}\quad m_q^{-1}\ l_q\quad m_q = 1$$

$S_g$

$$P_1^{-1}\quad Q_1^{-1}\quad P_1\ Q_1 \dots P_q^{-1}\quad Q_q^{-1}\ P_q\ Q_q$$

$\mu$

where $P_i, Q_i \in F_g'$, in two lines

$$l_1^{-1}\quad m_1^{-1}\quad l_1\ m_1 \dots l_g^{-1}\quad m_g^{-1}\ l_g\quad m_g = 1$$

$S_g$

$$c_1^{-1}\quad 1\quad c_1\quad 1\ \dots\ c_g^{-1}\quad 1\quad c_g\quad 1$$

$\mu_1$

$\gamma$-group $S_g$ with the exact isomorphism $\mu_1 \times \nu_1$ we will write in three lines.

$$l_1^{-1}\quad m_1^{-1}\quad l_1\ m_1 \dots l_g^{-1}\quad m_g^{-1}\ l_g\quad m_g = 1$$

$S_g$

$$c_1^{-1}\quad 1\quad c_1\quad 1\ \dots\ c_g^{-1}\quad 1\quad c_g\quad 1$$

$\mu_1$

$$L_1^{-1}\quad M_1^{-1}\quad L_1\ M_1 \dots L_g^{-1}\quad M_g^{-1}\ L_g\quad M_g$$

$\nu_1$

where components $c_1, \dots, c_g$ generate the group $F_g'$, and components $L_1, M_1, \dots, L_g, M_g$ generate the group $F_g''$.

Epimorphism $\mu_1$ determines (**1** is situated under $m_i^{\pm 1}$, $c_i^{\pm 1}$ is situated under $l_i^{\pm 1}$) $\gamma$-isomorphic of the group $S_g$.

With the epimorphism $\nu_1$ $\gamma$-group $S_g$ we will consider the *word*

$\nu_1 = L_1^{-1} M_1^{-1} L_1\ M_1, \dots, L_g^{-1} M_g^{-1} L_g M_g$ which is equal to **1** in $F_g''$.

In the word $\nu_1$ let's set any way of reduction letter by letter that converses the word $\nu_1$ in **1**, i.e. each pair of reducing letters $d_i^{\varepsilon}$, $d_i^{-\varepsilon}$ of the word $\nu_1$ we will connect with the *line* from above. The set of near located lines that connect self-inverse subwords of any components of the word $\nu_1$ we will call a *tape arch* or just *tape*. The subwords which are connected by a tape we will call *foundations* of this tape.

An arch is called *interior* if its foundations are self-inverse and belong to the same commutator.

An arch is called a *language* if its foundations belong to different commutators.

A commutator of $\gamma$-group $S_g$, $\mu_1 \times \nu_1$ in the form

$$\begin{array}{lllll} \dots l_i^{-1} & m_i^{-1} & l_i & m_i \dots & S_g \\ c_i^{-1} & 1 & c_i & 1 & \mu_1 \\ L_i^{-1} & M_i^{-1} & L_i & M_i & \nu_1 \end{array}$$

where components $\mathbf{L_i}$ , $\mathbf{M_i} \in \mathbf{F_g}''$ , is called an *exact* commutator.

A comutator $\gamma$- group $\mathbf{S_g}$, $\mu_N \times \nu_N$ in the form

$$\begin{array}{lllll} \ldots\ l_i^{-1} & m_i^{-1} & l_i & m_i \ldots & S_g \\ \quad c_i^{-1} & 1 & c_i & 1 & \mu_N \\ \quad L_i^{-1} & L_iM_iL_i^{-1} & L_i & M_i^{-1}r_i & \nu_N \end{array}$$

where $\mathbf{r}$ is the only language, and other components $\mathbf{L_i^{-1}}$, $\mathbf{L_iM_iL_i^{-1}}$, $\mathbf{L_i}$ ,$\mathbf{M_i^{-1}}$ are mutually reduced by interior arches, is called *normal*. In the normal commutator $\mathbf{r_i=M_iL_iM_i^{-1}L_i^{-1}}$ in $\mathbf{F_g}''$ .

-----------------------------------------------------------------------------------------------

* The proof was given by R. Grigorchuk and P. Kurchanov [3].

A commutator $\gamma$-group $\mathbf{S_g}$, $\mu_0 \times \nu_0$ in the form

$$\begin{array}{lllll} \ldots\ l_i^{-1} & m_i^{-1} & l_i & m_i \ldots & S_g \\ \quad c_i^{-1} & 1 & c_i & 1 & \mu_0 \\ \quad L_i^{-1} & M_i^{-1} & L_i & M_i & \nu_0 \end{array}$$

where for every $\mathbf{i=1,...,g}$ $\mathbf{L_iM_i = M_i L_i}$ in the group $\mathbf{F_g}''$ , is called *unite*.

**Lemma 2**. For $\gamma$-group $\mathbf{S_g}$ , $\mathbf{g>1}$ with exact epimorphism $\mu_1 \times \nu_1$ on group $\mathbf{F_g}' \times \mathbf{F_g}''$ such $\gamma$-epimorphism $\mathbf{l_1^{-1} \to l_i^{-1}m_1}$ & ...& $\mathbf{l_g^{-1} \to l_g^{-1}\, m_g}$ of the $\gamma$-group $\mathbf{S_g}$ on $\gamma$-group $\mathbf{S_g}$, such normal epimorphism $\mu_N \times \nu_N$ of the $\gamma$-group $\mathbf{S_g}$ on group $\mathbf{F_g}' \times \mathbf{F_g}''$ exist that diagram

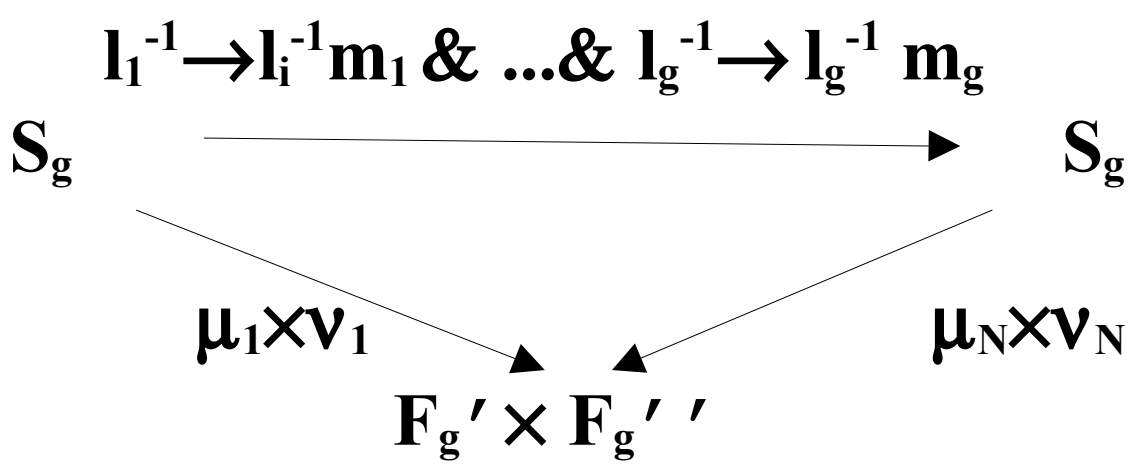

is commutative.

**Proof.** Let's consider $\gamma$-group $\mathbf{S_g}$ with the exact epimorphism $\mu_1 \times \nu_1$

$$\begin{array}{llllllll} l_1^{-1} & m_1^{-1} & l_1 & m_1 \ldots & l_g^{-1} & m_g^{-1} & l_g & m_g = 1 \\ c_1^{-1} & 1 & c_1 & 1 & c_g^{-1} & 1 & c_g & 1 \\ L_1^{-1} & M_1^{-1} & L_1 & M_1 & L_g^{-1} & M_g^{-1} & L_g & M_g \end{array}$$

Let's apply $\gamma$-isomorphism

$\mathbf{l_1^{-1} \to l_1^{-1}\, m_1}$ & …& $\mathbf{l_g^{-1} \to l_g^{-1}\, m_g}$

to $\gamma$-group $\mathbf{S_g}$ .

We will get

$$\begin{array}{llllllll} l_1^{-1} & m_1^{-1} & l_1 & m_1 \ldots & l_g^{-1} & m_g^{-1} & l_g & m_g = 1 \\ c_1^{-1} & 1 & c_1 & 1 & c_g^{-1} & 1 & c_g & 1 \\ L_1^{-1}M_1^{-1} & M_1^{-1} & M_1L_1 & M_1 & L_g^{-1}M_g^{-1} & M_g^{-1} & M_g L_g & M_g \end{array}$$

For every $\mathbf{i=1,...,g}$ using equalities $\mathbf{M_i\, L_i = A_i}$ , $\mathbf{M_i = A_i\, L_i^{-1}}$, $\mathbf{L_i^{-1} = B_i^{-1}A_i^{-1}}$ we construct components $\mathbf{A_i}$, $\mathbf{B_i}$ .

Next we will replace $\mathbf{M_i}$ to $\mathbf{A_iB_i^{-1}A_i^{-1}}$ & $\mathbf{L_i}$ to $\mathbf{A_iB_i}$ .

We will get $\gamma$-group $S_g$ with normal epimorphism

$$\begin{array}{cccccccccc}
l_1^{-1} & m_1^{-1} & l_1 & m_1 & \ldots & l_g^{-1} & m_g^{-1} & l_g & m_g =1 \\
c_1^{-1} & 1 & c_1 & 1 & & c_g^{-1} & 1 & c_g & 1 \\
A_1^{-1} & A_1B_1A_1^{-1} & A_1 & A_1B_1^{-1}A_1^{-1} & & A_g^{-1} & A_gB_gA_g^{-1} & A_g & A_gB_g^{-1}A_g^{-1}
\end{array}$$

**The composite epimorphisms for fragments $\gamma$-group $S_g$ with normal epimorphism.**

**1**. Let's consider a fragment that consists of two commutators **[1]** and **[i]** of $\gamma$-group $S_g$, with normal epimorphism $\mu_N \times \nu_N$

$$\begin{array}{cccccccccc}
l_1^{-1} & m_1^{-1} & l_1 & m_1 \cdot & l_i^{-1} & m_i^{-1} & l_i & m_i & \\
c_1^{-1} & 1 & c_1 & 1 & c_i^{-1} & 1 & c_i & 1 & [1]\cdot[i] \\
A_1^{-1} & A_1B_1A_1^{-1} & A_1 & B_1^{-1}r_1 & A_i^{-1} & A_iB_iA_i^{-1} & A_i & B_i^{-1}r_i &
\end{array}$$

where **i=2,...,g**. According to definition of commutator $\mu_N \times \nu_N$, the components $A_1^{-1}$, $B_1^{-1}$, $A_1$, $B_1$ of the commutator **[1]** and for **i=2,...,g** components $A_i^{-1}$, $B_i^{-1}$, $A_i$, $B_i$ of the commutator **[i]** completely reduce by interior arches of commutators. And languages $r_1$, $r_i$ can be determined from the equalities $r_1= B_1A_1B_1^{-1}A_1^{-1}$, $r_i= B_iA_iB_i^{-1}A_i^{-1}$.

With $\varphi(i)$ **(i=1,...,g)** we will denote epimorphism $m_1 \rightarrow m_1 (l_i^{-1}m_i^{-1} l_i m_i)^{-1}$.

Let's apply **epimorphism $\varphi(i)$** to the fragment **[1][i]** **(i=2,...,g)**.

We will get the fragment

$$\begin{array}{ccccccccc}
l_1^{-1} \cdot & l_i^{-1} & m_i^{-1} & l_i & m_i \cdot & m_1^{-1} & l_1 & m_1 \\
c_1^{-1} & c_i^{-1} & 1 & c_i & 1 & 1 & c_1 & 1 \\
A_1^{-1} & A_i^{-1} & A_iB_iA_i^{-1} & A_i & B_i^{-1}r_i & r_i^{-1}A_1B_1A_1^{-1} & A_1 & B_1^{-1}r_1r_i
\end{array}$$

where $r_i = B_iA_i B_i^{-1}A_i^{-1}$ ; **(i=1,...,g)**.

**2**. Let's consider the fragment that consists of the generator $l_1^{-1}$ and commutator **[i]** $\gamma$-group $S_g$ with normal epimorphism

$$\begin{array}{cccccc}
l_1^{-1} \cdot & l_i^{-1} & m_i^{-1} & l_i & m_i & \\
c_1^{-1} & c_i^{-1} & 1 & c_i & 1 & l_1^{-1}\cdot[i] \\
A_1^{-1} & A_i^{-1} & A_iB_iA_i^{-1} & A_i & B_i^{-1}r_i &
\end{array}$$

where **i=2,...,g**, $r_i = B_i A_iB_i^{-1}A_i^{-1}$

Let's define **epimorphism $\psi(i)$** **(i=2,...,g)**, consisted of four sequentially placed elementary isomorphisms

**1.i** $l_i^{-1} \rightarrow l_1 l_i^{-1}$ & $c_i^{-1} \rightarrow c_1 c_i^{-1}$

**2.i** $m_i \rightarrow l_1 m_i$

**3.i** $l_i \rightarrow l_1 l_i$ & $c_i \rightarrow c_1 c_i$

**4.i** $m_i^{-1}\to l_1\, m_i^{-1}$

Let's apply epimorphism $\psi(i)$ $(i=2,\dots,g)$ to the fragment $l_1^{-1}\times[i]$.

**1.i.** $l_i^{-1}\to l_1\, l_i^{-1}$ & $c_i^{-1}\to c_1\, c_i^{-1}$)

$$l_i^{-1}\qquad m_i^{-1}\qquad l_i\;\cdot\; l_1^{-1}\;\cdot\; m_i$$

$$c_i^{-1}\qquad 1\qquad c_i\qquad c_1^{-1}\qquad 1$$

$$A_1^{-1}A_i^{-1}\quad A_iB_iA_i^{-1}\quad A_i\,A_1\quad A_1^{-1}\quad B_i^{-1}r_i$$

**2.i.** $m_i\to l_1\, m_i$ (temporary breaks $\gamma$-isomorphic of the second line)

$$l_i^{-1}\qquad m_i^{-1}\;\cdot\; l_1^{-1}\;\cdot\; l_i\qquad m_i$$

$$c_i^{-1}\qquad c_1\qquad c_1^{-1}\qquad c_i\qquad c_1^{-1}$$

$$A_1^{-1}A_i^{-1}\quad A_iB_iA_i^{-1}A_1\quad A_1^{-1}\quad A_iA_1\quad A_1^{-1}B_i^{-1}r_i$$

**3.i.** $l_i\to l_1\, l_i$ & $c_i\to c_1\, c_i$

$$l_i^{-1}\;\cdot\; l_1^{-1}\;\cdot\; m_i^{-1}\qquad l_i\qquad m_i$$

$$c_i^{-1}\qquad c_1^{-1}\qquad c_1\qquad c_i\qquad c_1^{-1}$$

$$A_1^{-1}A_i^{-1}A_1\quad A_1^{-1}\quad A_iB_iA_i^{-1}A_1\quad A_1^{-1}A_iA_1\quad A_1^{-1}B_i^{-1}r_i$$

**4.i.** $m_i^{-1}\to l_1\, m_i^{-1}$ (returns $\gamma$-isomorphic to the second line)

We will get $\gamma$-isomorphic fragment

$$l_i^{-1}\qquad m_i^{-1}\qquad l_i\qquad m_i\;\cdot\; l_1^{-1}$$

$$c_i^{-1}\qquad 1\qquad c_i\qquad 1\qquad c_1^{-1}$$

$$A_1^{-1}A_i^{-1}A_1\quad A_1^{-1}A_iB_iA_i^{-1}A_1\quad A_1^{-1}A_iA_1\quad A_1^{-1}B_i^{-1}r_i\,A_1\quad A_1^{-1}$$

(Every $\psi(i)$ $(i=2,\dots,g)$ carries generator $l_1^{-1}$ of the fragment $l_1\times[i]$ to the right through the commutator $[i]$, and components $A_i$, $A_iB_iA_i^{-1}$, $B_i^{-1}r_i$ are joined by the component $A_1$ and convert respectively to $A_1^{-1}A_iA_1$, $A_1^{-1}A_iB_iA_i^{-1}A_1$, $A_1^{-1}B_i^{-1}r_i\,A_1$).

So epimorphism $\psi(i)$ $(i=2,\dots,g)$ is $\gamma$-isomorphism.

**<u>Lemma3</u>**. For $\gamma$-group $S_g$, $g>1$ with normal epimorphism $\mu_N\times\nu_N$ on group $F_g'\times F_g''$ we can indicate such $\gamma$-epimorphism $\rho$ of the $\gamma$-group $S_g$ on a free composition of groups $S_{g-1}*[l_1,m_1]=1$, such exact epimorphism $\mu_1\times\nu_1$ of the $\gamma$-group $S_{g-1}$ on $F_g'\times F_g''$, such unite epimorphism $\mu_0\times\nu_0$ of the group $[l_1,m_1]=1$ on $F_g'\times F_g''$ that diagram

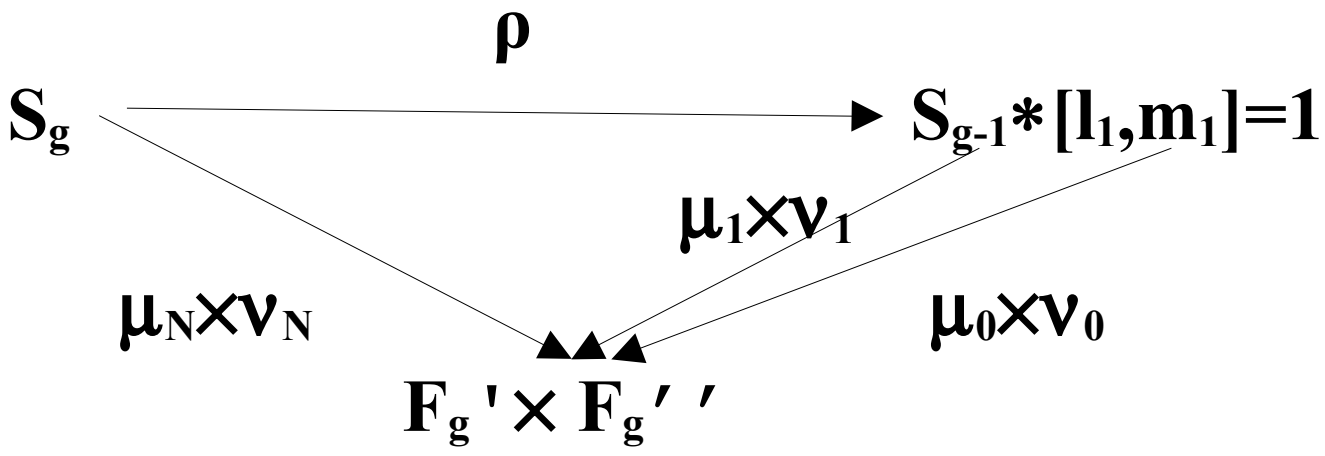

is commutative.

**Proof**.

Let's consider $\gamma$-group $S_g$ with normal epimorphism $\mu_N\times\nu_N$ on group $F_g'\times F_g''$

$$l_1^{-1}\quad m_1^{-1}\quad l_1\quad m_1\;\cdot\; l_2^{-1}\quad m_2^{-1}\quad l_2\quad m_2\cdot\ldots\cdot l_g^{-1}\quad m_g^{-1}\quad l_g\quad m_g=1$$

$$c_1^{-1} \quad 1 \quad c_1 \quad 1 \quad c_2^{-1} \quad 1 \quad c_2 \quad 1 \quad c_g^{-1} \quad 1 \quad c_g \quad 1$$
$$A_1^{-1}\ A_1B_1A_1^{-1}\ A_1\ \ B_1^{-1}r_1\ \ A_2^{-1}\ A_2B_2A_2^{-1}\ A_2\ \ B_2^{-1}r_2 \dots A_g^{-1}\ A_gB_gA_g^{-1}\ A_g\ \ B_g^{-1}r_g$$

where $r_1, r_2, \dots, r_g$ are languages. We will notice that because in every normal commutator **[i]** of the group $S_g$ there is an equality of foundations of interior arches, i.e. $A_i^{-1}\cdot A_iB_iA_i^{-1}A_iB_i^{-1}=1$ in the group $F_g'$, then the composition of languages $r_1 r_2 \dots r_g=1$ in the group $F_g''$.

Let's apply to $S_g$, $\mu_N\times\nu_N$ (consequently) composites of epimorphisms $\varphi(2),\dots,\varphi(g)$.

$$m_1 \rightarrow m_1 (l_2^{-1}m_2^{-1} l_2 m_2)^{-1}, \dots, m_1 \rightarrow m_1 (l_g^{-1}m_g^{-1} l_g m_g)^{-1}$$

(For every $i=2,\dots,g$ the transformation $\varphi(i)$ carries $m_1^{-1}l_1 m_1$ to the right through $l_i^{-1}m_i^{-1}l_i m_i$ ).

We will get $\gamma$-group **P**

$$l_1^{-1}\times l_2^{-1} \quad m_2^{-1} \quad l_2 \quad m_2\times\dots\times l_g^{-1} \quad m_g^{-1} \quad l_g \quad m_g \ \cdot \ m_1^{-1} \quad l_1 \quad m_1 =1$$
$$c_1^{-1}\ c_2^{-1} \quad 1 \quad c_2 \quad 1 \quad c_g^{-1} \quad 1 \quad c_g \quad 1 \quad 1 \quad c_1 \quad 1 \quad A_1^{-1}$$
$$A_2^{-1}\ A_2B_2A_2^{-1}\ A_2\ \ B_2^{-1}r_2 \dots A_g^{-1}\ \ A_gB_gA_g^{-1}\ A_g\ \ B_g^{-1}r_g\ \ A_1B_1A_1^{-1}\ A_1\ \ B_1^{-1}r_1r_2 \dots r_g$$

where $r_1, r_2, \dots, r_g$ are languages, $r_i = B_iA_i B_i^{-1}A_i^{-1}$ $(i=2,\dots,g)$.

Because in $S_g$, $\mu_N\times\nu_N$ $r_1 r_2 \dots r_g=1$, and in $\gamma$-isomorphisms $\varphi(2),\dots,\varphi(g)$ languages $r_1, r_2, \dots, r_g$ don't change, then in $\gamma$-group **P** $r_1 r_2 \dots r_g=1$ in $F_g''$.

To $\gamma$-group **P** (consequently) we will apply epimorphisms : $\psi(2)$ to the commutator **[2]** ,..., $\psi(g)$ to the commutator **[g].**

(Epimorphism $\psi(2)$&...& $\psi(g)$ consequently carries the generator **l1-1** through commutators **[2] ,..., [g]** . As a result **l1-1** with **m1-1 l1 m1** forms the commutator **[1]** . At the same time every component **V** of the third line $\gamma$-group **P** transfers to the component **A1-1VA1** ).

As a result we get $\gamma$-group in the form **[2]...[g][1]**

$$l_2^{-1} \quad m_2^{-1} \quad l_2 \quad m_2 \quad \dots \quad l_g^{-1} \quad m_g^{-1} \ = \ 1$$
$$c_2^{-1} \quad 1 \quad c_2 \quad 1 \quad c_g^{-1} \quad 1$$
$$A_1^{-1}A_2^{-1}A_1\ \ A_1^{-1}A_2B_2A_2^{-1}A_1\ \ A_1^{-1}A_2\ A_1\ \ A_1^{-1}B_2^{-1}r_2\ A_1 \dots A_1^{-1}\ A_g^{-1}A_1\ \ A_1^{-1}A_gB_gA_g^{-1}A_1$$
$$l_g \quad m_g \quad \times \quad l_1^{-1} \quad m_1^{-1} \quad l_1 \quad m_1$$
$$c_g^{-1} \quad 1 \quad c_1^{-1} \quad 1 \quad c_1 \quad 1$$
$$A_1^{-1}\ A_g\ A_1 \quad A_1^{-1}B_g^{-1}r_g\ A_1 \quad A_1^{-1}\ A_1B_1A_1^{-1} \quad A_1 \quad B_1^{-1}r_1r_2 \dots r_g$$

Because in $\gamma$- group **P** $r_1 r_2 \dots r_g=1$ in $F_g''$, and in $\gamma$-isomorphisms $\Psi(2),\dots,\Psi(g)$ languages $r_1, r_2, \dots, r_g$ aren't used, then in $\gamma$-group **[2]...[g][1]** we have $r_1r_2 \dots r_g=1$ in $F_g''$ .

We will apply to $\gamma$-group **[2]...[g][1]** epimorphism $r_1 r_2 \dots r_g \rightarrow 1$. Because in received $\gamma$-group in **[1]** components by generators $m_1^{-1}$ and $m_1$ are self-inverse then $A_1B_1A_1^{-1} = B_1$ in $F_g''$, and consequently, $A_1B_1 = B_1A_1$ in $F_g''$.

As a result the last commutator $\gamma$-group **[2]...[g][1]** will get the form

$$\begin{array}{ccccc} \ldots\; l_1^{-1} & m_1^{-1} & l_1 & m_1 \\ c_1^{-1} & 1 & c_1 & 1 \\ A_1^{-1} & B_1^{-1} & A_1 & B_1 \end{array} \qquad \mathbf{[1]}$$

where $A_1B_1 = B_1A_1$ in $F_g''$ and, consequently, **[1]** is a unite commutator.

In $\gamma$-group **[2]...[g]** for every $i=1,\ldots, g$ we will construct components $K_i$, $T_i$ such that

$$K_i = A_1^{-1}A_iA_1,\ T_i = A_1^{-1}A_i\, B_i^{-1}A_i^{-1}A_1 \quad (i=1,\ldots,g)$$

and make substitutions

$$A_1^{-1}A_iA_1 \rightarrow K_i \,\&\, A_1^{-1}A_i\, B_i^{-1}A_i^{-1}A_1 \rightarrow T_i$$

We will get $\gamma$-group $S_{g-1}$ in the form

$$\begin{array}{cccccccc} l_2^{-1} & m_2^{-1} & l_2 & m_2\ \ldots & l_g^{-1} & m_g^{-1} & l_g & m_g = 1 \\ c_2^{-1} & 1 & c_2 & 1 & c_g^{-1} & 1 & c_g & 1 \\ K_2^{-1} & T_2^{-1} & K_2 & T_2 & K_g^{-1} & T_g^{-1} & K_g & T_g \end{array}$$

Thus initial $\gamma$-group $S_g$ with normal epimorphism $\mu\times\nu$ by $\gamma$-isomorphisms transfers into a free composition of $\gamma$-group $S_{g-1}$ with exact epimorphism and commutator with a unite epimorphism

$$\begin{array}{cccccccccccc} l_2^{-1} & m_2^{-1} & l_2 & m_2\ \ldots & l_g^{-1} & m_g^{-1} & l_g & m_g = 1 * & l_1^{-1} & m_1^{-1} & l_1 & m_1 = 1 \\ c_2^{-1} & 1 & c_2 & 1 & c_g^{-1} & 1 & c_g & 1 & c_1^{-1} & 1 & c_1 & 1 \\ K_2^{-1} & T_2^{-1} & K_2 & T_2 & K_g^{-1} & T_g^{-1} & K_g & T_g & A_1^{-1} & B_1^{-1} & A_1 & B_1 \end{array}$$

where $A_1B_1 = B_1A_1$ in $F_g''$ .

**<u>Theoreme</u>**. For group $S_g$, $g>1$ with epimorphism $\mu\times\nu$ on group $F_g'\times F_g''$ we can indicate such epimorphism $\tau$ of the group $S_g$ on a free composition of groups $S_{g-1} * [l_1,m_1]=1$ , such an exact epimorphism $\mu_1\times\nu_1$ of the $\gamma$-group $S_{g-1}$ on $F_g'\times F_g''$, such unite epimorphism $\mu_0\times\nu_0$ of the group $[l_1,m_1]=1$ on group $F_g'\times F_g''$, that diagram

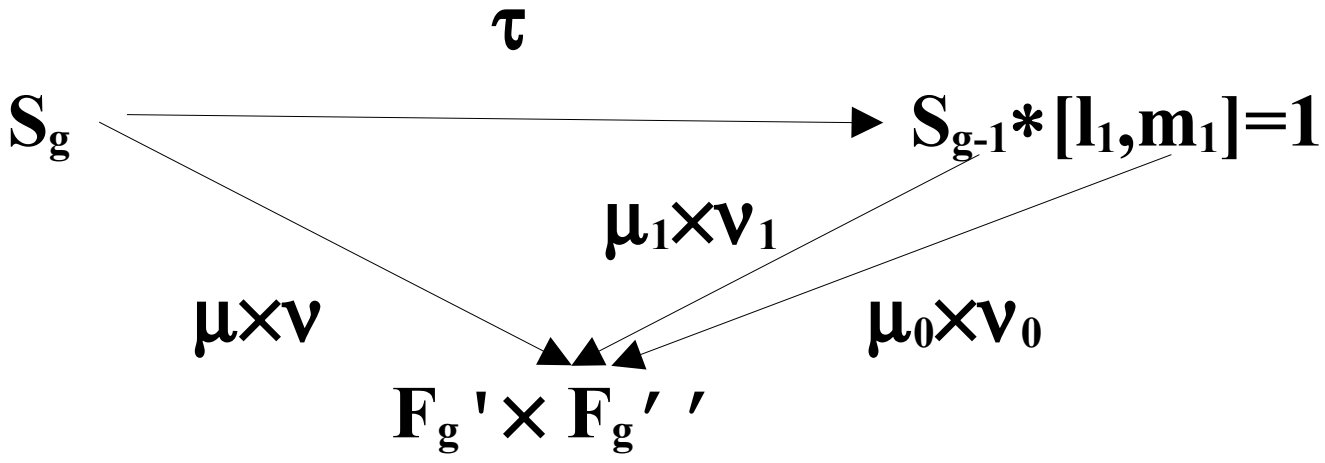

is commutative.

**<u>Proof</u>.** Combining consequently diagrams of the lemmas **1**,**2**,**3**, we will get that the diagram of the theorem is commutative.

## Literature.